\documentclass[11pt]{amsart}

\usepackage[english]{babel}
\usepackage{amsmath,amssymb, amsthm,amsbsy} 
\usepackage{hyperref}  
\hypersetup{colorlinks=true, linkcolor=blue, anchorcolor=blue, 
citecolor=red, filecolor=blue, menucolor=blue,
urlcolor=blue}

\numberwithin{equation}{section}
\newtheorem{theorem}{Theorem}[section]
\newtheorem{corollary}[theorem]{Corollary}

\theoremstyle{remark}
\newtheorem{remark}{Remark}[section]

\theoremstyle{definition}

\newcommand{\R}{{\mathbb R}}
 
\newcommand{\Z}{{\mathbb Z}}
\newcommand{\C}{{\mathbb C}}

\newcommand{\Hn}{\mathcal{H}_\nu}

\newcommand{\sola}{e^{itH^{a/2}}f}

\newcommand{\hank}{\mathcal{H}_\nu}
\newcommand{\four}{\mathcal{F}}

\newcommand{\heav}{\chi_{\mathbb{R}^+}}

\def\XXint#1#2#3{{\setbox0=\hbox{$#1{#2#3}{\int}$ }
\vcenter{\hbox{$#2#3$ }}\kern-.58\wd0}}

\begin{document}

\title
[Fractional Schr\"odinger with Aharonov-Bohm field]
{Weak dispersive estimates for fractional Aharonov-Bohm-Schr\"odinger groups}

\begin{abstract}
\end{abstract}

\date{\today}    


\author{F.~Cacciafesta}
\address{Federico Cacciafesta: Dipartimento di Matematica e Applicazioni, Universit$\grave{\text{a}}$ di Milano Bicocca, via R. Cozzi, 53, 20125 Milano, Italy.}
\email{federico.cacciafesta@unimib.it}

\author{L.~Fanelli}
\address{Luca Fanelli: SAPIENZA Universit$\grave{\text{a}}$ di Roma, Dipartimento di Matematica, P.le A. Moro 5, 00185 Roma, Italy.}
\email{fanelli@mat.uniroma1.it}


\thanks{
The two authors were supported by the Italian project FIRB 2012 {\it Dispersive Dynamics: Fourier Analysis and Calculus of Variations}.}

\begin{abstract}
We prove local smoothing, local energy decay and weighted Strichartz inequalities for fractional Schr\"odinger equations with a Aharonov-Bohm magnetic field in 2D. Explicit representations of the flows in terms of spherical expansions of the Hamiltonians are involved in the study. An improvement of the free estimate is proved, when the total flux of the magnetic field through the unit sphere is not an integer.
\end{abstract}

\date{\today}
\maketitle

\section{Introduction}

Given $\alpha\in\R$, consider the vector field (Aharonov-Bohm)
$$
A_B:\R^2\setminus\{(0,0)\}\to\R^2,
\qquad
A_B(x)=\left(-\frac{x_2}{|x|^2},\frac{x_1}{|x|^2}\right),
\qquad
x=(x_1,x_2)
$$
and the following quadratic form on $L^2(\R^2)$
$$
q_\alpha:\mathcal D(q_\alpha)\to[0,+\infty),
\qquad
q_\alpha[\psi]:=\int\left|\left(-i\nabla+A_B\right)\psi\right|^2\,dx,
$$
where the domain $\mathcal D(q_\alpha)$ is the completion of $\mathcal C^\infty_c(\R^2\setminus\{0\})$ with respect to the norm induced by $q_\alpha$. Since $q_\alpha$ is positive and symmetric, by the Friedrichs' Extension Theorem we can define the self-adjoint Hamiltonian 
\begin{equation}\label{eq:H}
H=\left(-i\nabla+\alpha\left(-\frac{x_2}{|x|^2},\frac{x_1}{|x|^2}\right)\right)^2,
\end{equation}
with its natural form domain, which coincides with the operator at the right-hand side of \eqref{eq:H} on $\mathcal C^\infty_c(\R^2\setminus\{0\})$. By Spectral Theorem, we can hence perform functional calculus on $H$ and define $f(H)$, being $f:\R\to\R$ any Borel-measurable function. In particular, we consider the positive powers $H^{a/2}$, $a>0$ and their associated Schr\"odinger flows $S(t):=e^{itH^{a/2}}$. The unitary (on $L^2$) group $S(t)$ uniquely defines the solution $u(t,\cdot):=e^{itH^{a/2}}f(\cdot)$ to the Cauchy problem
\begin{equation}\label{eq:schro}
\begin{cases}
\partial_tu = iH^{a/2}u
\\
u(0,\cdot)=f(\cdot)\in L^2(\R^2).
\end{cases}
\end{equation}
We will refer to \eqref{eq:schro} as the fractional Schr\"odinger equation with Aharonov-Bohm magnetic potential $A_B$. Notice that, when $a=2$, \eqref{eq:schro} is a Schr\"odinger equation with a magnetic potential $A_B$. On the other hand, when $a=1$, the flow $u(t,\cdot)=e^{itH^{1/2}}f(\cdot)$ has a clear connection with the solution of the following wave equation
$$
\begin{cases}
\partial_t^2v+Hv=0
\\
v(0,\cdot)=g(\cdot),
\\
\partial_tv(t,\cdot)=h(\cdot)
\end{cases}
$$
which is given by the formula
$$
v(t,\cdot)=\cos\left(tH^{\frac12}\right)g(\cdot)+\frac{\sin\left(tH^{\frac12}\right)}{H^{\frac12}}h(\cdot)
=\Re\left(e^{itH^{1/2}}\right)g(\cdot)+\frac{\Im\left(e^{itH^{1/2}}\right)}{H^{\frac12}}h(\cdot).
$$
Throughout this manuscript, we will refer to the case $\alpha=0$ as the {\it free case}. We remark that, as soon as $\alpha\in\mathbb Z$, $H_\alpha$ is unitarily equivalent to the free Hamiltonian $-\Delta$ (see e.g. \cite{pank-rich} and references therein). For this reason, from now on we will restrict to the case $\alpha\in[0,1)$.
 Among the many interesting features which the flow $e^{itH^a}$ enjoys, we first mention the invariance of equation \eqref{eq:schro} under the scaling $(t,x)\mapsto(\lambda^{-a}t,\lambda^{-1}x)$: for this reason, we can look at equation \eqref{eq:schro} as a critical (linear) perturbation of the free dispersive model. In recent years, critical perturbations of dispersive PDE's received a lot of interest, essentially motivated by the study of nonlinear models. The dispersive phenomenon can be quantified in several different ways in terms of a priori estimates for solutions. Time decay of $L^p$-norms has been recently studied and proved, for equation \eqref{eq:schro} with $a=2$ in \cite{fanfel1} and then generalized to a larger family of critical potentials in \cite{fanfel2, fanfel3, FGK, GK, K}.  As a consequence, Strichartz estimates can be obtained for \eqref{eq:schro} (with $a=2$) from the $L^1-L^\infty$ decay, by applying the standard Ginibre-Velo and Keel-Tao methods in \cite{GV, KT}. Nevertheless, neither sharp time decay estimates nor Strichartz estimates for \eqref{eq:schro} are known, at the best of our knowledge, when $a\neq2$. 
On the other hand, Strichartz estimates can be proved for critical 0-order perturbations of the free Schr\"odinger Hamiltonian without using the time decay, as shown in \cite{burqplanch, burqplanch2}, which are the crucial references of this manuscript.  
It is quite surprising that, in this case, perturbation techniques do hold. The strategy in \cite{burqplanch, burqplanch2} relies on a $TT^\star$-argument and a suitable mix of free Strichartz and local smoothing estimates \'a la Morawetz. It is easy to check that the argument fails in presence of a critical first-order potential as in \eqref{eq:schro}.
%

The aim of this paper is to prove local smoothing and local energy decay for solutions to \eqref{eq:schro}, in the same style as in \cite{burqplanch}. Before stating our main results, we briefly sketch a spectral picture of $H$ and introduce some notations (for more details we refer to \cite{adamtet,pank-rich}).
It is well known that $H$ is exactly solvable: its spectrum is purely absolutely continuous and coincides with the positive real axis, i.e. $\sigma(H)=\sigma_{ac}(H)=[0,+\infty)$. 
Throughout the paper, we will always use the canonical decomposition of $L^2(\R^2)$ in spherical harmonics. More precisely, given the complete orthonormal set on $L^2(\mathbb S^2)$ $\{\phi_m\}_{m\in\mathbb Z}$, with $\phi_m=\phi_m(\theta)=\frac{e^{im\theta}}{\sqrt{2\pi}}$, $\theta\in[0,2\pi)$,
 one has the canonical isomorphism
\begin{equation}\label{eq:isomorf}
L^2(\R^2)\cong \bigoplus_{m\in\Z}L^2(\mathbb{R}_+,rdr)\otimes [\phi_m]
\end{equation}
where we are denoting with $[\phi_m]$ the one dimensional space spanned by $\phi_m$ and with $\|f\|_{L^2_{rdr}}^2=\int_0^\infty|f(r)|^2rdr$. In this representation, the operator $H$ is equal to \cite[Sec. 2]{DS}
 \begin{equation}\label{dec}
H= \bigoplus_{m\in\mathbb{Z}}H_{\alpha,m}\otimes 1
 \end{equation}
 with 
 $$
 \displaystyle H_{\alpha,m}=-\frac{d^2}{dr^2}-\frac1r\frac{d}{dr}+\frac{(m+\alpha)^2}{r^2}.
 $$
 Therefore, the eigenvalue problem for $H$ leads to the Bessel equation, which can be solved after imposing boundary conditions, to obtain the generalized eigenfunctions
\begin{equation*}
\Psi_\alpha(r,\phi,k,\theta)=\sum_{m=-\infty}^{+\infty}i^{|m|}e^{im(\phi-\theta)}e^{i(\pi/2)(|m|-|m+\alpha|)}J_{|m+\alpha|}(kr).
\end{equation*}

We can now state the main result of this manuscript.
  \begin{theorem}\label{teo1}
Let $a>0$, $\alpha\in\mathbb R$ and $\varepsilon\in\left(0,\frac14+\frac12\text{dist}(\alpha,\mathbb Z)\right)$. Then for every $f\in L^2$ the following estimate holds
\begin{equation}\label{locsmooth}
\| |x|^{-\frac12-2\varepsilon} H^{\frac{a-1}4-\varepsilon}\sola\|_{L^2_tL^2_x}\leq C \|f\|_{L^2}
\end{equation}
with a constant $C$ depending on $\alpha$, $d$ and $\varepsilon$.

In addition, in the endpoint case $\varepsilon=0$ the following local estimate holds
\begin{equation}\label{locendec}
\sup_{R>0}R^{-1/2}\|\sola\|_{L^2_t L^2_{|x|\geq R}}\leq C\|H^\frac{1-a}4f\|_{L^2_x}.
\end{equation}
\end{theorem}

\begin{remark}
  Estimate \eqref{locsmooth} is false, also in the free case, for $\varepsilon=0$. On the other hand, it is interesting to notice that \eqref{locsmooth} fails for $\alpha\in\mathbb Z$ and $\varepsilon=\frac14$, in 2D. Indeed, the dimension $d=2$ is critical with respect to estimate \eqref{locsmooth}, with $\varepsilon=\frac14$, due to the fact that the weight $|x|^{-1}$ is too singular at the origin. Nevertheless, the presence of the field $A_B$, as it is well known, generically improves the angular ellipticity of $H$, if $\alpha\notin\mathbb Z$, and this usually permits to obtain better estimates than in the free case, as \eqref{locsmooth} shows. This phenomenon also appears for variational inequalities (see \cite{lw}), and for weighted dispersive estimates (see \cite{fanfel3, FGK, GK, K}).
  Roughly speaking, the higher the spherical frequency is, the better is the dispersive phenomenon we are measuring, as it will be clear in the sequel. The improvement arises since the introduction of the external potential is cutting the 0-frequency from the spectrum of the spherical operator.
 We finally remark that \eqref{locsmooth} holds with $\varepsilon=\frac14$, in dimensions larger or equal to 3 (see \cite{BVZ, burqplanch}).
\end{remark}

\begin{remark}
A frequency-dependent version of estimate \eqref{locsmooth} is inequality \eqref{restest} below. Indeed, we see that, after the decomposition \eqref{dec}, the restricted operator $H_{\alpha,m}$ satisfies \eqref{locsmooth} in the range $\varepsilon\in\left(0,\frac14+\frac12|m+\alpha|\right)$. Roughky speaking, the worst frequencies for \eqref{locsmooth} are the (at most) two closest integers to $\alpha$ (which in the free case coincide with the zero-frequency, as observed in \cite{burqplanch, burqplanch2} first).
\end{remark}

\begin{remark}
  We stress that inequalities \eqref{locsmooth} and \eqref{locendec} are the only dispersive estimates which are available for \eqref{eq:schro}, in the case $a\neq2$. Although we are still not able to use them to prove stronger dispersive inequalities, as Strichartz, they still represent a tool of independent interest. We also remark that the multiplier techniques of Morawetz type to prove local smoothing also fail in this case, since they usually requires the space dimension to be larger or equal to 3 (see e.g. \cite{BRV, DF, FV}). An analogous result has been recently proved for the Dirac equation with a Coulomb potential in \cite{cacser}.
\end{remark}

\begin{remark}
Theorem \ref{teo1} should be compared with Theorem 1.5 in \cite{fangwang}, in which the analogous estimates are proved for the the free flows (namely the case $\alpha=0$ in \eqref{eq:schro}). It should be noticed that in the free case an additional gain of angular derivative is shown; in fact, the same gain is expected to hold in our magnetic case as well, and this is suggested by the fact the constant defined in \eqref{constant} seems to allow some additional power of $\nu$,
and it will be subject of further investigation.
\end{remark}

\begin{remark}
  One may wonder if the distorted derivatives $H^{(1-a)/4}$ can be replaced by the usual derivatives $|D|^{(1-a)/4}$, both in estimate \eqref{locsmooth} and \eqref{locendec}. In dimension $d=2$ the usual Sobolev space $H^1(\R^2)$ is strictly bigger than the one generated by $H$, so that the answer is quite likely negative. 
\end{remark}

As a corollary, we can extend estimate \eqref{locsmooth} to the Klein-Gordon flow $e^{it\sqrt{H^a+1}}$.
\begin{corollary}\label{kg}
Let $a>0$, $\alpha\in\mathbb R$ and $\varepsilon\in\left(0,\frac14+\frac12\text{dist}(\alpha,\mathbb Z)\right)$. Then for every $f\in L^2$ the following estimate holds
\begin{equation}\label{locsmooth2}
\| |x|^{-\frac12-2\varepsilon} H^{\frac{a-1}4-\varepsilon}e^{it\sqrt{H^\frac{a}2+1}}\|_{L^2_tL^2_x}\leq C \|(H^\frac{a}2+1)^{\frac14}f\|_{L^2}
\end{equation}
with a constant $C$ depending on $\alpha$, $d$ and $\varepsilon$.
\end{corollary}
The proof immediately follows by Theorem \ref{teo1} and Theorem 2.4 in \cite{danckat}), and we omit further details. 

We conclude with a final application of Theorem \ref{teo1}.
Although we are still not able to prove Strichartz estimates for $e^{itH^{a/2}}$, it is still possible to prove some weighted version of them. In the following, we use the polar coordinates $x=r\omega$, $r\geq0$, $\omega\in\mathbb S^1$, and given a measurable function $F=F(t,x):\R\times\R^2\to\C$ we denote by
$$
\|F\|_{L^q_tL^q_{rdr}L^2_\omega}
:=
\int_{-\infty}^{+\infty}\left(r\int_{0}^{+\infty}\left(\int_{\mathbb S^1}|F(t,r,\omega)|^2\,d\sigma\right)^{q}\,dr\right)\,dt,
$$
being $d\sigma$ the surface measure on the sphere.
\begin{corollary}\label{weightstr}
Let $a>0$, $\alpha\in\mathbb R$, $\varepsilon\in\left(0,\frac14+\frac12\text{dist}(\alpha,\mathbb Z)\right)$ and $q\in[2,\infty]$. Then for every $f\in L^2$ the following estimate hold
\begin{equation}\label{wstr}
\|r^{-\frac12-2\varepsilon}e^{itH^{a/2}}f\|_{L^q_tL^q_{rdr}L^2_\omega}\leq
C
\| H^{\frac14+\varepsilon-\frac{a}{2q}}\Lambda_\omega^{\frac{1-b}2+\frac{b-1}r}f\|_{L^2},
\end{equation}
for some $C>0$, where $\Lambda_\omega=\sqrt{1-\Delta_\omega}$ and $\Delta_\omega$ is the Laplace-Beltrami operator on $S^1$ (we have introduced the polar coordinates $x=r\omega$).
\end{corollary}
The proof immediately follows by interpolation between estimate \eqref{locsmooth} and the 2D Sobolev inequality 
$$
\sup_{r>0} r^{\frac{1-4\varepsilon}2}\|f(r\omega)\|_{L^2_\omega}\leq C\| |D|^{\frac12+2\varepsilon}\Lambda_\omega^{-2\varepsilon}f\|_{L^2_x}
\leq C\| |H|^{\frac14+\varepsilon}\Lambda_\omega^{-2\varepsilon}f\|_{L^2_x},
$$
for $\varepsilon\in\left(0,\frac14+\frac\alpha2\right)$
together with the usual diamagnetic inequality.

The rest of the paper is devoted to the proof of Theorem \ref{teo1}.

\section{Proof of Theorem \eqref{teo1}}

Our proof follows closely the one of Theorems 1,2 in \cite{burqplanch}. 
In view of the decomposition \eqref{dec}, it is sufficient to obtain a suitable bound, uniform in $m$, for the projection $H_{\alpha,m}$. Let us fix $m\in\Z$ and denote by
$$
A_\nu:=H_{\alpha,m}=-\frac{d^2}{dr^2}-\frac1r\frac{d}{dr}+\frac{\nu^2}{r^2}
$$
where $\nu=|m+\alpha|>0$, which is a self-adjoint operator on the natural domain, for the same reasons as above (see \cite{pank-rich}).

Let us now introduce the standard Hankel transform of order $\nu>0$ as
$$
(\mathcal{H}_\nu\phi)(\xi)=\int_0^\infty J_{\nu}(r|\xi|)\phi(r\xi/|\xi|) r\:dr.
$$
The following properties are satisfied:
\begin{enumerate}
\item $\mathcal{H}_\nu^2=$Id;
\item $\mathcal{H}_\nu$ is self adjoint;
\item $\mathcal{H}_\nu$ is an $L^2$ isometry;
\item $\mathcal{H}_\nu A_\nu=|\xi|^2 \mathcal{H}_\nu$.
\end{enumerate}

For the proof of (1)-(2)-(3), see \cite{planchon}. 
The proof of (4) relies on the fact that the Bessel functions $J_\nu(r|\xi|)$ are (generalized) eigenfunctions for the restricted operator $A_\nu$. Indeed, given $f(x)=\psi(r)\phi_m(\theta)\in L^2(\mathbb{R}^+,rdr)\otimes[\phi_m]$, with $\phi_m(\theta)=(2\pi)^{-\frac12}e^{im\theta}$, the eigenvalue equation for $A_\nu$ reads
$$
\left(-\frac{d^2}{dr^2}-\frac1r\frac{d}{dr}+\frac{(m+\alpha)^2}{r^2}\right)\psi(r)=|\xi|^2\psi(r),
$$
and after the substitution $r\rightarrow |\xi|r$, we obtain the Bessel equation.
As a consequence, we have ($\langle,\cdot,\cdot\rangle$ denotes the standard $L^2$ product)
\begin{align*}
\mathcal{H}_\nu( A_\nu\phi)
&
=\langle J_\nu(r|\xi|), A_\nu\phi(r)\rangle
=\langle A_\nu J_\nu(r|\xi|),\phi\rangle
\\
&
=|\xi|\langle J_\nu(r|\xi|),\phi\rangle
\end{align*}
as $A_\nu$ is selfadjoint, and this proves (4).
Therefore, by (1)---(4)  we can write the fractional powers of $A_\nu$ as
\begin{equation}\label{fracpow}
A_\nu^{\sigma/2}= \mathcal{H}_\nu |\xi|^\sigma \mathcal{H}_\nu,
\end{equation}
or better as the following integral operators
$$
A^{\sigma/2}_\nu \phi(r,\theta)=\int_0^{\infty}k_{\nu,\nu}^\sigma \phi(s,\theta) s\:ds.
$$
Here the kernel $k_{\nu,\nu}$ is explicitly given by
\begin{equation}\label{kernel}
k_{\nu,\nu}^\sigma(r,s)=\begin{cases}
\displaystyle
\frac{2^{\sigma+1}\Gamma(\nu+\frac\sigma2+1)}{\Gamma(-\frac\sigma2)\Gamma(\nu+1)}
\frac{s^{\nu}}{r^{\sigma+\nu+2}}F(\nu+\frac\sigma2+1,\frac\sigma2+1;\nu+1;\frac{s^2}{r^2})\quad {\rm if} \: s<r;
\\
\displaystyle
\frac{2^{\sigma+1}\Gamma(\nu+\frac\sigma2+1)}{\Gamma(-\frac\sigma2)\Gamma(\nu+1)}
\frac{r^{\nu}}{s^{\sigma+\nu+2}}F(\nu+\frac\sigma2+1,\frac\sigma2+1;\nu+1;\frac{r^2}{s^2})\quad {\rm if}
\:r<s
\end{cases}
\end{equation}
(see \cite{planchon} for details).

We are now ready to prove Theorem \ref{teo1}. We claim that the following restricted version of estimate \eqref{locsmooth} holds
\begin{equation}\label{restest}
\| |x|^{-\frac12-2\varepsilon}A_\nu^{\frac{a-1}4-\varepsilon}S_\nu f\|_{L^2_tL^2_{rdr}}\leq C\|f\|_{L^2_{rdr}},
\end{equation}
with a constant $C>0$ independent on $\nu$ (hence on $m$)
where $S_\nu f=e^{itA_{\nu}^{a/2}}f$ is the unique solution of the initial value problem
\begin{equation*}\label{hankeq}
\begin{cases}
i\partial_t u+A_\nu^{a/2}  u=0\\
u(0,x)=f(x).
\end{cases}
\end{equation*}
By \eqref{fracpow} and the fact that $\mathcal H_\nu$ is an isometry on $L^2(\R^+;rdr)$, we can write
\begin{align}\label{isonew}
\| |x|^{-\frac12-2\varepsilon}A_\nu^{\frac{a-1}4-\varepsilon}S_\nu f\|_{L^2_tL^2_{rdr}}
&
=
\| \mathcal H_\nu|x|^{-\frac12-2\varepsilon}\mathcal H_\nu \mathcal H_\nu A_\nu^{\frac{a-1}4-\varepsilon}\mathcal H_\nu \mathcal H_\nu S_\nu f\|_{L^2_tL^2_{rdr}}
\\
&
=
\| A_\nu^{-\frac14-\varepsilon}|\xi|^{\frac{a-1}2-2\varepsilon}\mathcal{H}_\nu S_\nu f\|_{L^2_tL^2_{rdr}}
\nonumber.
\end{align}
As a consequence, by \eqref{isonew}, the claim \eqref{restest} is equivalent to the following
 \begin{equation}\label{eq:claim}
I:=\| A_\nu^{-\frac14-\varepsilon}|\xi|^{\frac{a-1}2-2\varepsilon}\mathcal{H}_\nu S_\nu f\|_{L^2_tL^2_{rdr}}\leq C\|\Hn f\|_{L^2_{rdr}}=C\| f\|_{L^2_{rdr}},
\end{equation}
for some constant $C>0$ independent on $\nu$.
The advantage is that $\Hn S_\nu f$ solves a much simpler system that is, due to \eqref{fracpow},
\begin{equation}\label{hankeq}
\begin{cases}
i\partial_t \Hn S_\nu f+ |\xi|^{a} \Hn S_\nu f =0\\
\Hn S_\nu f(0,\xi)=\Hn f(\xi)
\end{cases}
\end{equation}
which can be solved as
$$
\Hn S_\nu f(t,\xi)=e^{it|\xi|^{a}}(\Hn f)(\xi).
$$
Taking the time-Fourier transform and commuting, we see that 
$$
(\mathcal{F}_{t\rightarrow \tau} \Hn S_\nu f)(\tau,\xi)=(\Hn f)(\xi)\delta(\tau-|\xi|^{a}).
$$
Then we obtain, by Plancherel, that
\begin{align}\label{estpart}
&
I=
\left\|\int_0^\infty k_{\nu,\nu}^{-1/2-2\varepsilon}(|\xi|,s) \delta(\tau-s^{a})\Hn f(s\xi/|\xi) s^{1+\frac{a-1}2-2\varepsilon}\:ds\right\|_{L^2_\tau L^2_{\xi}}
\nonumber
\\
& \ \ \ = 
\left\| \frac12 \tau^{\frac1a(1+\frac{a-1}2-2\varepsilon)+\frac{1-a}a} k_{\nu,\nu}^{-1/2-2\varepsilon}(|\xi|,\tau^{1/a})(\Hn f)(\tau^{1/a}\xi/|\xi|) \right\|_{L^2_\tau L^2_{\xi}}.
\nonumber
\end{align}
Using polar coordinates in space and the change of variable $\omega=\tau^{1/a}$, we hence get
$$
I=
 \int_0^\infty \int_0^\infty\int_{S^1}\omega^{-4\varepsilon+2}( k_{\nu,\nu}^{-1/2-2\varepsilon})(\rho,\omega))^2 |(\Hn f)(\omega \theta)|^2 d\theta \rho \:d\rho d\omega.
$$
As now 
$$A_\nu^{-1/2-2\varepsilon}=A_\nu^{-1/4-\varepsilon}A_\nu^{-1/4-\varepsilon}$$ we have 
$$
k_{\nu,\nu}^{-1-4\varepsilon}(r,t)=\int_0^{\infty}k_{\nu,\nu}^{-1/2-2\varepsilon}(r,s)k_{\nu,\nu}^{-1/2-2\varepsilon}(s,t) \:sds
$$
and using this fact with the choice $r=t=\omega$ and $s=\rho$ we can write again
\begin{equation}\label{accc}
I=
\frac12\int_0^\infty \int_{S^1}\omega^{-4\varepsilon+2}k_{\nu,\nu}^{-1-4\varepsilon}(\omega,\omega)|(\Hn f)(\omega\theta)|^2 d\theta d\omega.
\end{equation}
By \eqref{kernel} (more precisely, the values of $k_{\nu,\nu}$ the diagonal $s=r$, that correspond to the values of a Gauss hypergeometric function in $z=1$) we obtain
\begin{align}
I
&
=\| A_\nu^{-\frac14-\varepsilon}|\xi|^{\frac{a-1}2-2\varepsilon}\mathcal{H}_\nu S_\nu f\|_{L^2_tL^2_{rdr}}=C_{\nu,\varepsilon}^2\int_0^\infty \int_{S^1} |(\Hn f)(\omega \theta)|^2  \omega\:d\omega d\theta
\\
&
=C_{\nu,\varepsilon}^2\left\|\mathcal H_\nu f\right\|_{L^2(rdr)}
=
C_{\nu,\varepsilon}^2\left\| f\right\|_{L^2(rdr)},
\nonumber
\end{align}
being the constant $C_{\nu,\varepsilon}$ given by
\begin{equation}\label{constant}
C_{\nu,\varepsilon}=2^{1/2-2\varepsilon}\sqrt{\pi \frac{\Gamma(\nu-2\varepsilon+1/2)\Gamma(4\varepsilon)}{\Gamma(\nu+2\varepsilon+1/2)\Gamma(2\varepsilon+1/2)^2}}.
\end{equation}
By \eqref{isonew}, to complete the proof of the claim \eqref{restest} we just need to check that $C_{}$
 is bounded with respect to the parameter $m$ (or equivalently $\nu$). To this aim, we notice that $C_{\nu,\varepsilon}$ is finite, provided $0<\varepsilon<1/4+\frac{\nu}2$ and it is a decreasing function of $\nu$. This completes the proof of \eqref{restest}. The proof of \eqref{locsmooth} now easily follows by \eqref{restest}, thanks to \eqref{eq:isomorf} and the $L^2$-orthogonality of the system $\{\phi_m\}_{m\in\mathbb Z}$, together with the fact that the weight $|x|^{-\frac12-2\varepsilon}$ is radial and that
 $$
 \min_{m\in\mathbb Z}\nu=\min_{m\in\mathbb Z}|m+\alpha|=\text{dist}(\alpha,\mathbb Z),
 $$

We now turn to the proof of \eqref{locendec}, which requires a slight modification of the above approach. By means of \eqref{eq:isomorf}, write $f\in L^2(\R^2)$ as $f=\sum_{m\in\mathbb Z}f_m(r,\theta)$, with $f_m(r,\theta)=\psi_m(r)\phi_m(\theta)$, $\psi_m\in L^2(\R^+;rdr)$. Setting $\xi=|\xi|\omega$, and using the change of variables $|\xi|^a=s$, we have, due to \eqref{hankeq},
\begin{align*}
e^{it H^{a/2}}f_m
&
= \hank  [e^{it |\xi|^a} \hank f_m]
\\
&
=
C\int_{S^1}\:\int_0^{+\infty}  e^{it|\xi|^a}J_\nu(r|\xi|) \hank f_m(|\xi|\omega) |\xi|\:d|\xi|\:d\sigma(\omega)
\\
&
=
C\int_{S^1}\:\int_0^{+\infty}e^{its} J_\nu(rs^{\frac1a}) \hank f_m(s^{\frac1a}\omega) s^{\frac2a-1}\:ds\:d\sigma(\omega)
\\
&
=
C\four_{s\rightarrow t} \left\{J_\nu(rs^{\frac1a}) \hank f_m(s^{\frac1a}\omega)s^{\frac2a-1}\heav\right\},
\end{align*}
with $C>0$ independent on $m$ and $f$,
being $\four_{s\rightarrow t}$ the Fourier transform in the $s$-variable.
We now take the $L^2_t L^2_{|x|\leq R}$ norm and apply Plancherel, to get
\begin{align}
&
\|e^{it H^a}f_m\|_{L^2_t L^2_{|x|\leq R}}^2
=
C\left\|\four_{s\rightarrow t} \left\{J_\nu(rs^{\frac1a}) \hank f_m(s^{\frac1a}\omega)s^{\frac2a-1}\heav\right\}\right\|_{L^2_t L^2_{rdr(0,R)}}^2
\label{finalled}
\\
& \ \ \ 
=
C\left\|J_\nu(rs^{\frac1a}) \hank f_m(s^{\frac1a}\omega)s^{\frac2a-1}\heav\right\|_{L^2_sL^2_{rdr(0,R)}}^2
\nonumber
\\
& \ \ \ 
=
C\int_{S^1}\int_0^{+\infty}\left( \int_0^{R}  J_\nu(rs^{\frac1a})^2r dr\right) |\hank f_m(s^{\frac1a}\omega)|^2 s^{2(\frac2a-1)} dsd\sigma(\omega)
\nonumber
\\
& \ \ \ 
=
C\int_{S^1}\int_0^{+\infty}\left( \int_0^{R}  J_\nu(r|\xi|)^2r dr\right) |\hank f_m(|\xi|\omega)|^2 |\xi|^{1+2-a} d|\xi|d\sigma(\omega).
\nonumber
\end{align}
We now notice that, for every $R>0$,
\begin{equation}\label{estbesint}
\int_0^{R}  J_\nu(r|\xi|)^2r  dr\leq \frac{CR}{|\xi|}  
\end{equation}
for some constant $C>0$ independent on $\nu$ (see \cite[pag. 63]{strich}). As a consequence, by \eqref{finalled} we obtain 
\begin{align}
\||e^{it H^a}f_m\|_{L^2_t L^2_{|x|\leq R}}^2
&
\leq CR \int_{S^1}\int_0^{+\infty}|\hank f_m(|\xi|\omega)|^2 |\xi|^{1+(1-a)} d|\xi|d
\label{finalled2}
\\
&=
 CR\|A_\nu^{\frac{1-a}4} f_m\|_{L^2},
 \nonumber
\end{align}
with a constant $C>0$ independent on $\nu$.
Estimate \eqref{locendec} now easily follows by \eqref{finalled2} and the decomposition $f=\sum_{m\in\mathbb Z}f_m=\sum_{m\in\mathbb Z}\psi_m(r)\phi_m(\theta)$, together with the $L^2$-orthogonality of the set $\phi_m$. The proof of Theorem \ref{teo1} is now complete.


\begin{thebibliography}{9}

\bibitem{adamtet}
R. Adami and A. Teta.
\newblock On the Aharonov-Bohm Hamiltonian.
\newblock {\em Letters in Math. Phys.} 43, 45-54 (1998).

\bibitem{BRV}
J. Barcel\'o, A. Ruiz, and L. Vega.
\newblock Some dispersive estimates for Schr\"odinger equations with repulsive potentials.
\newblock {\em J. Func. Anal.} 236 (2006), 1-24.

\bibitem{BVZ}
J. Barcel\'o, L. Vega, and M. Zubeldia.
\newblock The forward problem for the electromagnetic Helmholtz equation with critical singularities.
\newblock {\em Adv. Math.} 240 (2013), 636-671.

\bibitem{burqplanch}
N. Burq, F. Planchon, J.G. Stalker and A.S. Tahvildar-Zadeh. 
\newblock Strichartz estimates for the wave and Schr\"odinger equations with the inverse-square potential. 
\newblock {\em J. Funct. Anal.} 203 (2), 519--549 (2003). 

\bibitem{burqplanch2}
N. Burq, F. Planchon, J.G. Stalker and A.S. Tahvildar-Zadeh. 
\newblock Strichartz estimates for the Wave and Schr\"odinger Equations with Potentials of Critical Decay.
\newblock {\em Indiana Univ. Math. J.} 53 (6) (2004), 1665-1680.

\bibitem{cacser}
F. Cacciafesta and Eric S\'er\'e.
\newblock Local smoothing estimates for the Dirac Coulomb equation in 2 and 3 dimensions.
\newblock {http://arxiv.org/abs/1503.00945}

\bibitem{DS}
L. Dabrowski, and P. Stovicek.
\newblock Aharonov-Bohm effect with $\delta$-type interaction.
\newblock {\em J. Math. Phys.} 39 (1998), 47-62.

\bibitem{danckat}
P. D'Ancona.
\newblock Kato smoothing and Strichartz estimates for wave equations with magnetic potentials.
\newblock {\em Comm. Math. Phys}

\bibitem{DF}
P. D'Ancona, and L. Fanelli.
\newblock Smoothing estimates for the Schr\"odinger equation with unbounded potentials.
\newblock {\em Journ. Diff. Eq.} 246 (2009), 4552-4567.

\bibitem{fanfel1}
L. Fanelli, V. Felli, M. Fontelos and A. Primo.
\newblock Time decay of scaling critical electromagnetic Schr\"odinger flows,
\newblock {\em Comm. Math. Phys.} 324 (3), 1033--1067 (2013). 

\bibitem{fanfel2}
L. Fanelli, V. Felli, M. Fontelos and A. Primo.
\newblock Time decay of scaling invariant electromagnetic Schrödinger equations on the plane.
\newblock {\em Comm. Math. Phys. } 337 (2015), 1515-1533.

\bibitem{fanfel3}
L. Fanelli, V. Felli, M. Fontelos and A. Primo.
\newblock Frequency-dependent time decay of Schr\"odinger flows.
\newblock To appear in {\em J. Spectral Theory}.

\bibitem{FGK}
L. Fanelli, G. Grillo, and H. Kovarik.
\newblock Improved time-decay for a class of scaling-critical Schr\"odinger flows.
\newblock {\em J. Func. Anal.} 269 (2015), 3336-3346.

\bibitem{FV} 
L. Fanelli, and L. Vega.
\newblock Magnetic virial identities, weak dispersion and Strichartz inequalities.
\newblock {\em Math. Ann.} 344 (2009), 249-278.

\bibitem{fangwang}
D. Fang, and C. Wang.
\newblock Weighted Strichartz estimates with angular regularity and their applications.
\newblock {\em Forum Math.} 23 (2011), no. 1, 181-205.

\bibitem{GV}
J. Ginibre, and G. Velo.
\newblock Generalized Strichartz
inequalities for the wave equation. 
\newblock {\em J. Funct. Anal.} 133 no. 1
(1995), 50-68.

\bibitem{GK}
G. Grillo, and H. Kovarik.
\newblock Weighted dispersive estimates for two-dimensional Schr\"odinger operators with Aharonov-Bohm magnetic field.
\newblock {\em Journ. Diff. Eq.} 256 (2014), 3889-3911.

\bibitem{KT}
M. Keel, and T. Tao.
\newblock Endpoint Strichartz estimates.
\newblock {\em Am. J.
Math.} 120 no. 5 (1998), 955-980.

\bibitem{K}
H. Kovarik.
\newblock Resolvent expansion and time decay of the wave functions for two-dimensional magnetic Schr\"odinger operators.
\newblock {\em Comm. Math. Phys}. 337 (2015), 681-726.

\bibitem{lw}
A. Laptev, and T. Weidl
\newblock Hardy inequalities for
      magnetic Dirichlet forms.
      \newblock {\em Mathematical results in quantum
    mechanics (Prague, 1998)}, 299-305; \textit{Oper. Theory Adv. Appl.} 108,
    Birkh\"auser, Basel, 1999.

\bibitem{pank-rich}
K. Pankrashkin, and S. Richard.
\newblock Spectral and scattering theory for the Aharonov-Bohm operators. 
\newblock {\em Rev. Math. Phys.} 23 (2011), no. 1, 53?81.

\bibitem{planchon}
F. Planchon, J. Stalker, and A. S. Tahvildar-Zadeh.
\newblock  $L^p$ estimates for the wave equation with the inverse-square potential. 
\newblock {\em Discrete Contin. Dynam. Systems,} 9(2):427-442, (2003).

\bibitem{strich}
R. Strichartz.
\newblock Harmonic analysis as spectral theory of the Laplacians.
\newblock {\em J. Func. Anal.} 87 (1989), 51-148.



\end{thebibliography}
\end{document}